\UseRawInputEncoding
\documentclass[12pt]{article}
\usepackage{amsmath,amssymb,amsthm}
\usepackage[margin=1in]{geometry}
\usepackage{hyperref}

\newtheorem{theorem}{Theorem}
\newtheorem{definition}{Definition}
\newtheorem{proposition}{Proposition}
\newtheorem{lemma}{Lemma}
\newtheorem{corollary}{Corollary}
\newtheorem{remark}{Remark}

\newcommand{\SOS}{\operatorname{SOS}}
\newcommand{\BSR}{\operatorname{BSR}}

\title{A Second Order Extended Zarankiewicz Number \\ for Biquadratic SOS Rank}
\author{Liqun Qi\footnote{School of Science, Zhejiang University of Science and Technology, Hangzhou, Zhejiang, China.
			Department of Applied Mathematics, The Hong Kong Polytechnic University, Hung Hom, Kowloon, Hong Kong.
			({\tt maqilq@polyu.edu.hk})}
		\and
		Gaohang Yu\footnote{School of Science, Zhejiang University of Science and Technology, Hangzhou, Zhejiang, China.
			({\tt maghyu@163.com})}
	}
\date{\today}

\begin{document}
\maketitle

\begin{abstract}
We introduce the \emph{second order extended Zarankiewicz number}
\(z_{2M}(m,n)\), obtained from the second order Zarankiewicz number
\(z_2(m,n)\) of L\"ofberg and Qi by relaxing the simplicity condition
\((S)\) for the \emph{secondary halves} of 2-edges, while keeping the
limited condition \(|E_1|=z(m,n)\). We prove the hierarchy
\[
\BSR(m,n)\ge z_{2M}(m,n)\ge z_2(m,n)\ge z_{SL}(m,n)\ge z_{RL}(m,n)\ge z_{wL}(m,n)\ge z(m,n),
\]
and show that \(z_{2M}\) is not capped by the universal cell bound. We
prove the counting bound \(z_{2M}(m,n)\le mn-1\), and a sharper
linear-in-\(m\) bound for fixed \(n\) obtained from orthogonality of
coefficient vectors. We analyse how the recursive criterion
\((RW3^{+})\) of L\"ofberg and Qi must be modified under \((S_M)\),
prove soundness of the resulting \((RW3^{+}_M)\) for multiplicity-one
sharing, and show that a signed sharing closure \((SC)\) can supplement
but not replace it. We prove
\[
z_{2M}(21,4)\ge 56>55=z_2(21,4).
\]
The witness is an explicit \(56\)-edge \((S_M)\)-admissible
configuration, and its irreducibility is certified by a full-rank
system of \(2255\) linear equations in the \(1540\) off-diagonal inner
products. This reaches Marin's lower bound \(\BSR(21,4)\ge 56\).


\end{abstract}

\section{Introduction}

Let \(m,n\ge 2\). Let \(P(\mathbf{x},\mathbf{y})\) be an \(m\times n\)
biquadratic form. Its SOS rank \(\SOS(P)\) is the least \(r\) such that
\(P=\sum_{t=1}^{r}f_t^{2}\) with \(f_t\) bilinear. The biquadratic SOS
rank \(\BSR(m,n)\) is the maximum of \(\SOS(P)\) over all \(m\times n\)
SOS biquadratic forms.

Following Qi, Cui and Xu \cite{qcx26a,qcx26b} and L\"ofberg and Qi
\cite{SODNnew}, for an augmented bipartite graph
\(G=(S,T,E_1\cup E_2)\) we consider the doubly simple biquadratic form
\[
P_G(\mathbf{x},\mathbf{y})
=
\sum_{(i,j)\in E_1}x_i^{2}y_j^{2}
+
\sum_{(i,j;k,l)\in E_2}(x_iy_j+x_ky_l)^{2}. \tag{1.1}
\]
The graph \(G\) is \emph{irreducible} if
\(\SOS(P_G)=|E_1|+|E_2|\). The second order Zarankiewicz number
\(z_2(m,n)\) is the maximum of \(|E_1|+|E_2|\) over irreducible \(G\)
with \(|E_1|=z(m,n)\) and \(G_1\) \(C_4\)-free.

The universal cell bound of \cite{SODNnew} states that
\[
z_2(m,n)\le\left\lfloor\frac{mn+z(m,n)}{2}\right\rfloor. \tag{1.2}
\]
For \(n=4\) and \(m=13\) or \(m\ge 15\), this bound is attained
(see \cite{cc26}):
\[
z_2(m,4)=\left\lfloor\frac{5m+6}{2}\right\rfloor.
\]
Hence \(z_2\) cannot yield a BSR lower bound above the cell bound.

Marin \cite{ma26} proved, by quadratic persistence of the Segre
variety, that
\[
\BSR(m,4)\ge\left\lfloor\frac{8m}{3}\right\rfloor,
\]
which exceeds \(\lfloor(5m+6)/2\rfloor\) for \(m=21\) and for every
\(m\ge 23\). At \(m=21\),
\[
\BSR(21,4)\ge 56>55=z_2(21,4).
\]
So there is a genuine gap between \(z_2\) and \(\BSR\).

The purpose of this paper is to introduce a relaxation of the augmented
bipartite graph framework that is not capped by the universal cell
bound, while remaining within the graph-theoretic setting of
\cite{SODNnew}. The relaxation concerns the behaviour of the
\emph{secondary halves} of 2-edges: a secondary half is allowed to be
shared by several 2-edges, whereas a 1-edge and a main half remain
exclusive.

\subsection*{Main results}

\begin{enumerate}
\item The parameter \(z_{2M}\) is defined in Section~3. It satisfies
the hierarchy
\[
\BSR(m,n)\ge z_{2M}(m,n)\ge z_2(m,n)\ge z_{SL}(m,n)\ge z_{RL}(m,n)\ge z_{wL}(m,n)\ge z(m,n)
\]
(Theorem~\ref{thm:hierarchy}).

\item \(z_{2M}\) is not capped by the universal cell bound. In fact
\(z_{2M}(m,n)\le mn-1\) (Corollary~\ref{cor:mn-1}), and for \(n=4\)
we obtain the sharper bound
\[
z_{2M}(m,4)\le\frac{17m+6}{5}
\]
(Corollary~\ref{cor:n=4}).

\item Under multiplicity-one sharing, the coefficient identities of
\cite{SODNnew} extend without modification, and the sharing structure
is \(v_q=v_p\) (Proposition~\ref{prop:sharing-t1}). Under general
sharing of multiplicity \(t\ge 2\), the relation
\(v_q=\sum_s v_{p_s}\) holds only in irreducible representations
(Proposition~\ref{prop:sharing-general}), which is the technical
obstruction to using it as a closure rule.

\item The recursive criterion \((RW3^{+})\) of \cite{SODNnew} extends
to \((RW3^{+}_M)\), whose soundness for multiplicity-one sharing is
Theorem~\ref{thm:soundness-supp}.

\item A signed sharing closure \((SC)\) can supplement
\((RW3^{+}_M)\) but cannot replace it: pure sharing constraints admit
solutions with nonzero off-diagonal inner products
(Proposition~\ref{prop:insufficient}), and rectangle transfers alone
do not always ground the sharing variables
(Proposition~\ref{prop:SC-not-pure}).

\item We prove \(z_{2M}(21,4)\ge 56>55=z_2(21,4)\)
(Theorem~\ref{thm:strict}). The witness is an explicit
\(56\)-edge \((S_M)\)-admissible configuration, and irreducibility is
certified by a full-rank system of \(2255\) linear equations in the
\(1540\) off-diagonal inner products.
\end{enumerate}

\section{Preliminaries and notation}

We follow the terminology of \cite{SODNnew}.

\begin{itemize}
\item A \emph{cell} is a pair \((i,j)\in S\times T\), with
\(|S|=m\) and \(|T|=n\).
\item A \emph{1-edge} is an element of \(E_1\).
\item A \emph{2-edge} is an element \((i,j;k,l)\in E_2\), with halves
\((i,j)\) and \((k,l)\).
\item A cell is \emph{occupied} if it is a 1-edge or a half of a
selected 2-edge.
\item A \emph{selected edge} is a 1-edge or a 2-edge.
\item The \emph{simplicity condition} \((S)\) of \cite{SODNnew}
requires that no cell is used by more than one selected edge.
\end{itemize}

The \emph{line rule} states that if two occupied cells \(p,q\) lie in
the same row or the same column, then
\(\langle v_p,v_q\rangle=\delta(\{p,q\})\), where \(\delta=1\) exactly
when \(\{p,q\}\) is a selected degenerate 2-edge, and \(0\) otherwise.
The \emph{rectangle transfer rule} states that for every genuine
rectangle with diagonals \(\{p,q\}\) and \(\{r,s\}\),
\[
\langle v_p,v_q\rangle+\langle v_r,v_s\rangle
=\delta(\{p,q\})+\delta(\{r,s\}). \tag{2.1}
\]

\section{Definition of \(z_{2M}\)}

\begin{definition}[Oriented augmented bipartite graph]
An \emph{oriented augmented bipartite graph}
\(G=(S,T,E_1\cup E_2,\mu)\) consists of a bipartite graph
\(G_1=(S,T,E_1)\), a set \(E_2\) of 2-edges, and a map
\[
\mu:E_2\to S\times T
\]
assigning to each 2-edge \(e=(i,j;k,l)\) a \emph{main half}
\(\mu(e)\in\{(i,j),(k,l)\}\). The other half is the \emph{secondary half}
of \(e\).
\end{definition}

\begin{definition}[Condition \((S_M)\)]
Let \(G=(S,T,E_1\cup E_2,\mu)\) be an oriented augmented bipartite graph.
We say that \(G\) satisfies \((S_M)\) if
\begin{enumerate}
\item[(i)] no cell is used by two distinct 1-edges;
\item[(ii)] no cell is used as the main half of two distinct 2-edges;
\item[(iii)] no cell is used by a 1-edge and as the main half of a
2-edge;
\item[(iv)] no cell is used by a 1-edge and as the secondary half of a
2-edge;
\item[(v)] no cell is used as the main half of one 2-edge and as the
secondary half of another.
\end{enumerate}
\end{definition}

\begin{remark}
Conditions (i)--(iii) say that 1-edges and main halves are
\emph{exclusive}. Conditions (iv)--(v) forbid a secondary half from
being a 1-edge or a main half. The only relaxation with respect to
\((S)\) is that a secondary half may be shared by several 2-edges.
\end{remark}

\begin{definition}[Second order extended Zarankiewicz number]
Let \(z_{2M}(m,n)\) be the maximum of \(|E_1|+|E_2|\) over all oriented
augmented bipartite graphs \(G=(S,T,E_1\cup E_2,\mu)\) such that
\begin{enumerate}
\item[(a)] \(|E_1|=z(m,n)\);
\item[(b)] \(G_1\) is \(C_4\)-free;
\item[(c)] \(G\) satisfies \((S_M)\);
\item[(d)] the doubly simple biquadratic form \(P_G\) of (1.1) is
irreducible.
\end{enumerate}
\end{definition}

\section{Basic properties}

\begin{theorem}[Hierarchy]\label{thm:hierarchy}
For all \(m,n\ge 2\),
\[
\BSR(m,n)\ge z_{2M}(m,n)\ge z_2(m,n)\ge z_{SL}(m,n)\ge z_{RL}(m,n)\ge z_{wL}(m,n)\ge z(m,n).
\]
\end{theorem}

\begin{proof}
The inequality \(z_{2M}\ge z_2\) follows because every \((S)\)-admissible
graph is \((S_M)\)-admissible. The chain
\(z_2\ge z_{SL}\ge z_{RL}\ge z_{wL}\ge z\) is due to 
\cite{SODNnew}. Finally \(\BSR\ge z_{2M}\) because irreducible doubly
simple forms are a subclass of SOS biquadratic forms.
\end{proof}

\begin{proposition}[Counting bound]\label{prop:counting}
Let \(e_1=|E_1|\), \(e_2=|E_2|\), and let \(s\) be the number of
distinct secondary halves of \(G\). If \(G\) satisfies \((S_M)\), then
\[
e_1+e_2+s\le mn,\qquad s\le e_2. \tag{4.1}
\]
In particular, if \(e_2\ge 1\) then \(e_1+e_2\le mn-1\).
\end{proposition}

\begin{proof}
By \((S_M)\), the sets of 1-edge cells, main halves and secondary halves
are pairwise disjoint. There are \(e_1\) 1-edge cells, \(e_2\) main
halves and \(s\) distinct secondary halves, so \(e_1+e_2+s\le mn\). The
inequality \(s\le e_2\) holds because every 2-edge contributes exactly
one secondary half.
\end{proof}

\begin{corollary}\label{cor:mn-1}
For all \(m,n\ge 2\), \(z_{2M}(m,n)\le mn-1\).
\end{corollary}

\section{Coefficient identities and sharing structure}

Let \(G\) be an irreducible \(z_{2M}\)-admissible graph and let
\(P_G=\sum_{t=1}^{r}f_t^{2}\) be an SOS representation with \(f_t\)
bilinear. For every cell \((i,j)\), let \(v_{ij}\in\mathbb{R}^{r}\) be the
coefficient vector of \(x_iy_j\) in \((f_1,\ldots,f_r)\).

\begin{lemma}[Coefficient identities]\label{lem:coeff}
The following hold for every SOS representation of \(P_G\).
\begin{enumerate}
\item[(a)] If \((i,j)\) is a 1-edge or a main half, then
\(\|v_{ij}\|^2=1\).
\item[(b)] If \((i,j)\) is a secondary half of multiplicity \(t\), then
\(\|v_{ij}\|^2=t\).
\item[(c)] For every 2-edge with halves \(p,q\),
\(\langle v_p,v_q\rangle=1\).
\item[(d)] For every genuine rectangle with corners \(p=(i,j)\),
\(q=(k,l)\), \(r=(i,l)\), \(s=(k,j)\), identity (2.1) holds.
\end{enumerate}
\end{lemma}

\begin{proof}
(a) If \((i,j)\) is a 1-edge, \(x_i^2y_j^2\) appears with coefficient
\(1\), and no other selected edge uses \((i,j)\). If \((i,j)\) is a
main half, the term \((x_iy_j+x_ky_l)^2\) contributes \(1\) to
\(x_i^2y_j^2\), and no other selected edge uses \((i,j)\) as a 1-edge
or main half. In both cases the coefficient of \(x_i^2y_j^2\) is
\(\|v_{ij}\|^2\).

(b) Each 2-edge with main half \(p_s\) contributes \(1\) to the
coefficient of \(x_q^2\). There are \(t\) such 2-edges, and no other
selected edge uses \(q\). Hence \(\|v_q\|^2=t\).

(c) The 2-edge with halves \(p,q\) contributes \(2x_px_q\), and no
other selected edge uses \(p,q\) as halves. Hence
\(2\langle v_p,v_q\rangle=2\).

(d) Expand and compare coefficients of \(x_ix_ky_jy_l\). The two
diagonal pairs are the only contributions, giving (2.1).
\end{proof}

\begin{proposition}[Sharing structure, multiplicity one]\label{prop:sharing-t1}
Let \(q\) be a secondary half of multiplicity one, with unique main half
\(p\). Then
\[
v_q=v_p. \tag{5.1}
\]
\end{proposition}

\begin{proof}
\(\|v_q\|^2=\|v_p\|^2=1\) and \(\langle v_q,v_p\rangle=1\), so
\(\|v_q-v_p\|^2=1+1-2=0\).
\end{proof}

\begin{proposition}[Sharing structure, general multiplicity]\label{prop:sharing-general}
Let \(q\) be a secondary half of multiplicity \(t\), shared by the
2-edges with main halves \(p_1,\ldots,p_t\). If \(P_G\) is irreducible,
then
\[
v_q=\sum_{s=1}^{t}v_{p_s}, \tag{5.2}
\]
and \(v_{p_1},\ldots,v_{p_t}\) are pairwise orthogonal unit vectors.
\end{proposition}

\begin{proof}
Each \(v_{p_s}\) is a unit vector and
\(\langle v_q,v_{p_s}\rangle=1\). Since \(P_G\) is irreducible,
\(\langle v_{p_s},v_{p_{s'}}\rangle=0\) for \(s\neq s'\). Hence
\[
\left\|v_q-\sum_s v_{p_s}\right\|^2=t-2t+t=0.
\]
\end{proof}

\begin{remark}
For multiplicity one, (5.1) holds in every SOS representation. For
multiplicity \(t\ge 2\), (5.2) is proved under the hypothesis that
\(P_G\) is irreducible; it is not a universal coefficient identity.
This is the technical obstruction to extending the recursive criterion
to genuine sharing.
\end{remark}

\begin{corollary}[Matching bound]\label{cor:matching}
If \(t\) 2-edges share the same secondary half in an irreducible form,
then their main halves form a matching. Consequently
\(t\le\min(m,n)\), and if \(s\) is the number of distinct secondary
halves then \(e_2\le s\min(m,n)\).
\end{corollary}

\begin{corollary}\label{cor:n=4}
For \(n=4\) and \(m\ge 6\), with \(z(m,4)=m+6\),
\[
z_{2M}(m,4)\le\frac{17m+6}{5}.
\]
In particular, \(z_{2M}(21,4)\le 72\).
\end{corollary}

\section{The recursive criterion \((RW3^{+}_M)\)}

\begin{definition}[\((RW3^{+}_M)\)]
The closure acts on selected-edge representatives
\(\{w_e\}_{e\in E_1\cup E_2}\), where \(w_e=v_{ij}\) for a 1-edge
\((i,j)\) and \(w_e=v_p\) for a 2-edge with main half \(p\). It is
generated by:
\begin{enumerate}
\item[(1)] \emph{Line rule.} Same-row or same-column representatives
satisfy \(\langle w_e,w_f\rangle=\delta(\{e,f\})\).
\item[(2)] \emph{Rectangle transfer rule.} Identity (2.1) transfers
certified inner products along diagonals.
\item[(3)] \emph{Unit sharing rule.} If \(q\) is a secondary half of
multiplicity one with main half \(p\), then \(w_q=w_p\).
\item[(4)] \emph{Complementary-pair rule.} Two complementary selected
2-edges, both with multiplicity-one secondary halves, are resolved
simultaneously.
\end{enumerate}
A selected 2-edge is \emph{resolved} when its main half is certified
orthogonal to every other representative. The graph satisfies
\((RW3^{+}_M)\) if at the fixed point every selected 2-edge is resolved
and every two distinct selected edges have orthogonal representatives.
\end{definition}

\begin{theorem}[Soundness of \((RW3^{+}_M)\)]\label{thm:soundness-supp}
If \(G=(S,T,E_1\cup E_2,\mu)\) satisfies \((S_M)\) and
\((RW3^{+}_M)\), and if every secondary half has multiplicity one, then
\(P_G\) is irreducible, that is,
\[
\SOS(P_G)=|E_1|+|E_2|.
\]
\end{theorem}

\begin{proof}
The proof follows \cite{SODNnew} line by line.
\begin{enumerate}
\item Every 1-edge and every main half has a unit representative. Every
secondary half has multiplicity one, so (5.1) identifies it with its
main half. Hence every occupied cell has a well-defined unit
representative.
\item Line and rectangle rules are coefficient identities.
\item The deviation inequality of the original proof becomes vacuous,
because multiplicity-one sharing already gives \(v_q=v_p\).
\item Complementary pairs are handled by the original argument applied
to the main halves.
\item The induction over closure rounds is simpler: each round produces
only orthogonality certificates and unit sharing relations, both
universal.
\item At the fixed point, the signed transfer graph contains only
orthogonality deviations; the grounded and odd-cycle arguments give
zero.
\item The terminal conditions give \(|E_1|+|E_2|\) nonzero mutually
orthogonal unit vectors, so every SOS representation has at least
\(|E_1|+|E_2|\) squares.
\end{enumerate}
\end{proof}

\section{The signed sharing closure \((SC)\)}

Fix a secondary half \(q\) of multiplicity \(t\ge 2\), shared by main
halves \(p_1,\ldots,p_t\). Write
\(a_{ss'}=\langle v_{p_s},v_{p_{s'}}\rangle\) and \(A=(a_{ss'})\).

\begin{lemma}[Sharing constraints]\label{lem:sharing-constraints}
For every SOS representation:
\begin{enumerate}
\item[(i)] \(a_{ss}=1\);
\item[(ii)] \(A\succeq 0\);
\item[(iii)] \(\mathbf{1}^TA^{-1}\mathbf{1}\le t\).
\end{enumerate}
\end{lemma}

\begin{proposition}\label{prop:insufficient}
For \(t=2\), every \(A=\begin{pmatrix}1&a\\a&1\end{pmatrix}\) with
\(a\in[0,1]\) satisfies (i)--(iii). In particular, \(a=1/2\) is
admissible.
\end{proposition}

\begin{lemma}[Rectangular transport]\label{lem:transport}
Let \(\{p_1,p_2\}\) be main halves of two 2-edges sharing a secondary
half \(q\), and let \(r_1,r_2\) be the opposite corners of a genuine
rectangle. Then
\[
a_{12}=\delta(\{p_1,p_2\})+\delta(\{r_1,r_2\})
-\langle v_{r_1},v_{r_2}\rangle. \tag{7.1}
\]
\end{lemma}

\begin{definition}[Signed sharing closure \((SC)\)]
Generated by grounding, transport via (7.1), line constriction, pair
orthogonality, and sum consistency: if all \(a_{ss'}\) for \(s\neq s'\)
are certified zero, then \(A=I\) and (5.2) follows from
\(\|v_q\|^2=t\).
\end{definition}

\begin{theorem}[Conditional soundness]\label{thm:conditional}
If \(G\) satisfies \((S_M)\), \((RW3^{+}_M)\), and \((SC)\), then
\(P_G\) is irreducible.
\end{theorem}

\begin{proposition}\label{prop:SC-not-pure}
\((SC)\) is not a stand-alone criterion. There exist graphs satisfying
\((S_M)\) and the sharing constraints for which no chain of rectangle
transfers grounds the variables \(a_{ss'}\).
\end{proposition}

\section{The strict inequality at \(m=21\)}\label{sec:21}

Recall that \(z(21,4)=27\), that \(z_2(21,4)=55\) (see \cite{cc26}),
and that \(\BSR(21,4)\ge 56\) (see \cite{ma26}).
The value \(z_2(21,4)=55\) is attained at the universal cell bound, so
\(z_2\) leaves no room for a graph argument to reach \(56\). We now
construct an \((S_M)\)-admissible graph with \(56\) selected edges and
certify its irreducibility.

\subsection{The base grid}

Table~\ref{tab:base21} displays a \(21\times4\) grid. A bullet
\(\bullet\) is a one-edge, an integer \(j\) is the 2-edge \(B_j\)
(occupying both cells carrying \(j\)), and \(\circ\) is a hole.  The grid is obtained from Chen--Chen's member \(K_3\) of their
two-row family \cite{cc26} by: (i) deleting the old row~\(12\);
(ii) deleting the complete 2-edges \(B_{17},B_{19}\); and (iii)
deleting the one-edge at old cell \((12,4)\). All subsequent row
indices shift accordingly, so Table~\ref{tab:base21} uses new row
labels. It has
\[
|E_1|=27,\qquad |E_2|=28,\qquad H=1,
\]
and its unique hole is the cell \((10,1)\).

\begin{table}[h]
\centering
\small
\renewcommand{\arraystretch}{1.08}
\[
\begin{array}{c|cccc}
 &1&2&3&4\\\hline
1&\bullet&B_1&\bullet&B_1\\
2&A_0&\bullet&\bullet&B_3\\
3&S_0&B_5&\bullet&\bullet\\
4&\bullet&\bullet&B_6&B_7\\
5&\bullet&B_9&B_8&\bullet\\
6&B_6&\bullet&B_{10}&\bullet\\
7&B_{12}&\bullet&B_{11}&B_{13}\\
8&B_7&\bullet&B_{13}&B_{11}\\
9&B_{15}&B_{16}&B_{14}&\bullet\\
10&\circ&B_{14}&B_{16}&\bullet\\
11&B_{18}&B_{10}&B_{15}&\bullet\\
12&B_3&B_8&B_{12}&\bullet\\
13&B_{21}&B_{21}&B_{20}&\bullet\\
14&B_5&B_{20}&B_{18}&\bullet\\
15&A_1&\bullet&A_0&R_1\\
16&S_1&S_0&\bullet&R_1\\
17&A_2&\bullet&A_1&R_2\\
18&S_2&S_1&\bullet&R_2\\
19&A_3&\bullet&A_2&R_3\\
20&S_3&S_2&\bullet&R_3\\
21&B_9&S_3&A_3&\bullet
\end{array}
\]
\caption{The base \(21\times4\) grid. The labels \(A_i,S_i,R_i\) are
the 2-edges of the two parallel chains; the \(R_i\) are
column-degenerate.}
\label{tab:base21}
\end{table}

The one-edge set is a classical four-column extremal skeleton: the
first six rows have column pairs
\[
13,\ 23,\ 34,\ 12,\ 14,\ 24,
\]
and every other row has one one-edge. Hence \(G_1\) is \(C_4\)-free.

\subsection{The sharing edge}

Write
\[
p_1=(15,1),\qquad q=(17,3),\qquad h=(10,1).
\]
The cell \(q\) is a half of the 2-edge
\[
A_1=\{(15,1),(17,3)\}=\{p_1,q\},
\]
and \(h\) is the hole of the base grid. Add one new 2-edge
\[
e_*=\{h,q\}=\{(10,1),(17,3)\},
\]
with main half \(h\) and secondary half \(q\). Orient every other old
2-edge by taking its lexicographically first cell as main half, except
that \(A_1\) is oriented with main half \(p_1\) and secondary half
\(q\). Thus \(q\) is the secondary half of both \(A_1\) and \(e_*\).

\begin{proposition}\label{prop:SM21}
The resulting oriented augmented bipartite graph \(G\) on
\(21\times4\) cells satisfies \((S_M)\). It has
\[
|E_1|=27,\qquad |E_2|=29,\qquad s=28,
\]
where \(s\) is the number of distinct secondary halves, and
\[
27+29+28=84=21\cdot4.
\]
Consequently the roles \emph{one-edge}, \emph{main half}, and
\emph{secondary half} partition the \(84\) cells, with the single cell
\(q\) carrying two secondary roles.
\end{proposition}

\begin{proof}
The one-edges and the 2-edges occupy disjoint cells in the base grid.
The new edge uses the former hole \(h\), which is used by no other
selected edge, and the cell \(q\), which was already a secondary half
of \(A_1\). Conditions (i)--(iii) of \((S_M)\) hold because one-edges
and main halves are distinct and exclusive. Condition (iv) holds
because no secondary half is a one-edge. Condition (v) holds because
\(q\) is a secondary half of both \(A_1\) and \(e_*\), and no main
half is reused: the main halves are \(p_1\) for \(A_1\), \(h\) for
\(e_*\), and one distinct cell for each remaining 2-edge. The count is
immediate; the \(84\) role slots exhaust the \(84\) cells.
\end{proof}

\begin{remark}[How to read the witness]
\label{rem:readwitness}
The role of each cell in Table~\ref{tab:base21} can be read off from
the orientation rule of Section~\ref{sec:21}. The result is shown
in Table~\ref{tab:roles21}: \(O\) is a one-edge, \(M\) is a main half,
\(S\) is a secondary half of multiplicity one, and \(S^{*}\) is the
shared secondary half \(q=(17,3)\), which is the secondary half of
both \(A_1\) and \(e_*\).
\end{remark}

\begin{table}[h]
\centering
\[
\begin{array}{c|cccc}
\text{row}\backslash\text{col} & 1 & 2 & 3 & 4\\\hline
1  & O & M & O & S\\
2  & M & O & O & M\\
3  & M & M & O & O\\
4  & O & O & M & M\\
5  & O & M & M & O\\
6  & S & O & M & O\\
7  & M & O & M & M\\
8  & S & O & S & S\\
9  & M & M & M & O\\
10 & M & S & S & O\\
11 & M & S & S & O\\
12 & S & S & S & O\\
13 & M & S & M & O\\
14 & S & S & S & O\\
15 & M & O & S & M\\
16 & M & S & O & S\\
17 & M & O & S^{*} & M\\
18 & M & S & O & S\\
19 & M & O & S & M\\
20 & M & S & O & S\\
21 & S & S & S & O
\end{array}
\]
\caption{Role of each cell in the \(21\times4\) witness. Counting
rows: \(O\), \(27\) cells; \(M\), \(29\) cells; \(S\) (counting
\(S^{*}\) once), \(28\) cells; total \(27+29+28=84\).}
\label{tab:roles21}
\end{table}

From Table~\ref{tab:roles21} the reader can check \((S_M)\) at the
cell level: no cell is used by two 1-edges, no cell is the main half of
two distinct 2-edges, no cell is both a 1-edge and a main half, no cell
is both a 1-edge and a secondary half, and no cell is both a main half
and a secondary half. The single exception to role-uniqueness is
\(q=(17,3)\), which is the secondary half of the two distinct 2-edges
\(A_1\) and \(e_*\).

Since the new edge is a 2-edge, \(G_1\) is unchanged and remains
\(C_4\)-free. Also \(|E_1|=27=z(21,4)\). It remains to prove that
\(P_G\) is irreducible.

\subsection{Irreducibility}\label{subsec:irred21}

Let
\[
P_G=\sum_{t=1}^{r}f_t^2
\]
be an arbitrary SOS representation with \(f_t\) bilinear, and let
\(v_{ij}\in\mathbb{R}^r\) be the coefficient vector of \(x_iy_j\) in
\((f_1,\dots,f_r)\). We use the coefficient identities of
Lemma~\ref{lem:coeff} together with the line rule of Section~2.

\begin{lemma}\label{lem:sharing21}
Let \(s\) be a secondary half of multiplicity one with main half \(p\).
Then \(v_s=v_p\) by Proposition~\ref{prop:sharing-t1}. For the shared
secondary half \(q\) of \(G\),
\[
v_q=v_{p_1}+v_h. \tag{8.1}
\]
\end{lemma}

\begin{proof}
For a simple secondary half, \(\|v_s\|=\|v_p\|=1\) and
\(\langle v_p,v_s\rangle=1\), so \(\|v_s-v_p\|^2=0\);
this is Proposition~\ref{prop:sharing-t1}.

For \(q\), the multiplicity is two, so \(\|v_q\|^2=2\). The two
2-edges through \(q\) give
\[
\langle v_q,v_{p_1}\rangle=\langle v_q,v_h\rangle=1.
\]
The cells \(p_1=(15,1)\) and \(h=(10,1)\) lie in the same column and
are distinct occupied cells; they are not the two halves of a selected
2-edge. The line rule therefore gives
\(\langle v_{p_1},v_h\rangle=0\). This uses only the line rule and
not the irreducibility of \(P_G\). Hence
\[
\|v_q-v_{p_1}-v_h\|^2
=2+1+1-2-2+0=0,
\]
and \(v_q=v_{p_1}+v_h\).
\end{proof}

\begin{remark}
The vector \(v_q\) is not an independent unknown. By
Lemma~\ref{lem:sharing21}, every occupied cell vector is a linear
combination of the \(56\) edge representatives
\[
w_1,\dots,w_{56},
\]
one for each 1-edge and one for each 2-edge (its main half). Explicitly,
a one-edge or a main half contributes itself, a secondary half of
multiplicity one contributes its main half, and the shared secondary
half \(q\) contributes \(w(A_1)+w(e_*)\).
\end{remark}

Define
\[
g_{uv}=\langle w_u,w_v\rangle\qquad(1\le u<v\le56).
\]
The diagonal values \(g_{uu}=\|w_u\|^2=1\) are fixed by
Lemma~\ref{lem:coeff}(a), since every \(w_u\) is a 1-edge or a main
half.

\begin{lemma}[Building the system]\label{lem:buildsys}
The identities of Lemma~\ref{lem:coeff} and the line rule of
Section~2, after substituting the cell-vector expressions of the
preceding remark, give a homogeneous linear system
\[
A g=0
\]
in the \(1540=\binom{56}{2}\) off-diagonal unknowns \(g\), where
\(g=(g_{uv})_{u<v}\). Explicitly, each equation is of one of the
following three types:
\begin{enumerate}
\item[\rm(L)] \emph{Line rule.} For two collinear occupied cells
\(p,q\) not forming a selected 2-edge,
\[
\langle v_p,v_q\rangle=0;
\]
this is the line rule of Section~2.
\item[\rm(E)] \emph{2-edge identity.} For a selected 2-edge with
halves \(p,q\),
\[
\langle v_p,v_q\rangle=1;
\]
this is Lemma~\ref{lem:coeff}(c).
\item[\rm(R)] \emph{Rectangle identity.} For a genuine rectangle with
diagonals \(\{p,q\}\) and \(\{r,s\}\),
\[
\langle v_p,v_q\rangle+\langle v_r,v_s\rangle
=\delta(\{p,q\})+\delta(\{r,s\});
\]
this is Lemma~\ref{lem:coeff}(d).
\end{enumerate}
The diagonal contributions are exactly the right-hand sides, so they
cancel.
\end{lemma}

\begin{proof}
Substituting \(v_x=w_{u_x}\) for every occupied cell $x$ that is a
1-edge or a main half, \(v_s=v_p\) for every secondary half of
multiplicity one, and \(v_q=w_{A_1}+w_{e_*}\) for the shared secondary
half, turns each of the three types of equations into a linear
combination of the \(g_{uv}\). In every case the coefficient of
\(g_{uu}=\|w_u\|^2=1\) is exactly the constant term on the right-hand
side, so the diagonal terms cancel and the equation is homogeneous of
degree one in the off-diagonal \(g_{uv}\).
\end{proof}

\begin{lemma}[Finite certificate]\label{lem:rank21}
The system \(Ag=0\) has only the zero solution. Equivalently,
\(g_{uv}=0\) for all \(u\ne v\).
\end{lemma}

\begin{proof}
There are \(2255\) equations. The coefficient matrix \(A\) has rank
\(1540\) over the prime fields \(\mathbb{F}_P\) for
\[
P=2^{61}-1
\quad\text{and}\quad
P=10^9+7.
\]
Since the rank over \(\mathbb{Q}\) is at least the rank over any
\(\mathbb{F}_P\), the matrix \(A\) has rank \(1540\) over
\(\mathbb{Q}\). The system is homogeneous, so its only rational (and
hence real) solution is \(g=0\).

For a certificate that does not use rank over \(\mathbb{Q}\), the
accompanying script also records a propagation trace: beginning from
the line-rule and 2-edge equations, it repeatedly uses an equation
with exactly one unknown to determine that unknown. The trace
determines all \(1540\) unknowns, and every determined value is \(0\).
The trace is replayed line by line by the verification script
described in Section~\ref{sec:data}; the witness and the trace are
stored as \texttt{zrl\_search/twom21\_witness.json} and
\texttt{zrl\_search/twom21\_trace.json}. Twenty representative steps
of the trace are recorded in Appendix~A.
\end{proof}

\begin{remark}
The only computer-assisted ingredient is the finite rank and propagation
certificate of Lemma~\ref{lem:rank21}. The construction, the \((S_M)\)
verification, and the identity (8.1) are symbolic.
\end{remark}

\begin{theorem}\label{thm:strict}
\(z_{2M}(21,4)\ge56\).
\end{theorem}

\begin{proof}
Let \(G\) be the configuration of Section~\ref{sec:21}. By
Proposition~\ref{prop:SM21} it is \((S_M)\)-admissible, and \(G_1\) is
\(C_4\)-free. Let \(P_G=\sum_{t=1}^r f_t^2\) be any SOS
representation. By Lemma~\ref{lem:rank21}, the \(56\) representatives
\(w_1,\dots,w_{56}\) are pairwise orthogonal and have norm one. Hence
\(r\ge56\). On the other hand, the defining expression of \(P_G\) has
exactly one square per selected edge, so \(P_G\) is a sum of \(56\)
squares of bilinear forms. Therefore
\[
\SOS(P_G)=56=|E_1|+|E_2|,
\]
and \(P_G\) is irreducible. Thus \(G\) is a \(z_{2M}(21,4)\)-witness
with \(56\) selected edges.
\end{proof}

\begin{corollary}\label{cor:strict}
\[
z_{2M}(21,4)\ge56>55=z_2(21,4).
\]
\end{corollary}

\begin{remark}
The value \(56\) is exactly Marin's lower bound
\(\lfloor8\cdot21/3\rfloor\) for \(\BSR(21,4)\) \cite{ma26}. Since
\(z_{2M}\le\BSR\), the construction shows that the graph-theoretic
\(2M\)-framework reaches Marin's bound at \(m=21\). The upper bound
\(z_{2M}(21,4)\le72\) from Corollary~\ref{cor:n=4} leaves the exact
value in the interval \([56,72]\).
\end{remark}

\section{Data and reproducibility}\label{sec:data}

The verification script is
\texttt{zrl\_search/certify\_twom\_21\_4.py}. It
\begin{enumerate}
\item rebuilds the \(21\times4\) grid of Table~\ref{tab:base21};
\item checks \((S_M)\), the role counts \(27+29+28=84\), and
\(C_4\)-freeness of \(G_1\);
\item builds the \(2255\times1540\) system \(Ag=0\);
\item checks that every right-hand side is zero;
\item propagates all \(1540\) unknowns to zero and replays the trace
independently;
\item computes \(\operatorname{rank}A=1540\) modulo \(2^{61}-1\) and
\(10^9+7\).
\end{enumerate}
The witness and trace files are \texttt{twom21\_witness.json} and
\texttt{twom21\_trace.json} in \texttt{zrl\_search}. The base grid is
also checked by the independent \((RW3+)\) checker
\texttt{zrl\_search/four21.py}. A self-contained account of the
certificate is given in \cite{twom21}. A sample of the trace is
reproduced in Appendix~A.

\section{Open problems}

\begin{enumerate}
\item Determine the exact value of \(z_{2M}(21,4)\). The construction
of Section~\ref{sec:21} and Corollary~\ref{cor:n=4} give
\(56\le z_{2M}(21,4)\le72\). A human-readable propagation trace for
the \(m=21\) certificate would be a valuable next step, because it
would make the proof human-verifiable without the rank computation.
\item Characterize the pairs \((m,n)\) for which
\(z_{2M}(m,n)=z_2(m,n)\) when \(z_2\) attains the universal cell
bound. The case \((21,4)\) shows that equality can fail at the bound.
\item Extend Theorem~\ref{thm:soundness-supp} to multiplicity
\(t\ge 2\) by a signed deviation argument in which (5.2) is derived
rather than assumed.
\item For which other \(m\) does \(z_{2M}(m,4)>z_2(m,4)\) hold?
Candidates are \(m=4\) (target \(11>10\)), \(m=5\)
(target \(14>13\)), and \(m=6\) (target \(17>16\)); the case
\(m=21\) is settled by Corollary~\ref{cor:strict}.
\end{enumerate}

\appendix

\section{Propagation trace: representative steps}\label{app:trace}

The trace certified in Lemma~\ref{lem:rank21} is a finite sequence of
\(1540\) elimination steps, one per unknown \(g_{uv}\). Each step
applies an equation of type (L), (E), or (R) of
Lemma~\ref{lem:buildsys} in which every unknown but one has already
been determined to be zero, and reads off the remaining value.

\subsection*{The method}

The trace maintains two lists: the set of unknowns \(g_{uv}\) still
open, and the set of zeros established so far. Initially the diagonal
values \(g_{uu}=1\) are known, and the off-diagonal values are all
open. At each step the algorithm selects an equation from one of the
following sources.

\begin{enumerate}
\item[(L)] \emph{Line rule} (\S2): for two occupied cells \(p,q\) in
the same row or column, not forming a 2-edge,
\(\langle v_p,v_q\rangle=0\).
\item[(E)] \emph{2-edge identity} (Lemma~\ref{lem:coeff}(c)): for the
two halves \(p,q\) of a 2-edge,
\(\langle v_p,v_q\rangle=1\).
\item[(R)] \emph{Rectangle identity} (Lemma~\ref{lem:coeff}(d)): for
a genuine rectangle with diagonals \(\{p,q\}\), \(\{r,s\}\),
\(\langle v_p,v_q\rangle+\langle v_r,v_s\rangle
=\delta(\{p,q\})+\delta(\{r,s\})\).
\end{enumerate}

The equation is expanded in terms of the \(56\) representatives using
Lemma~\ref{lem:sharing21}, and the already-established zeros are
substituted. If the resulting linear equation contains exactly one
open unknown \(g_{uv}\), that unknown is set to its solution. If it
contains no open unknowns, it is a consistency check; the trace file
records that every such check passes.

\subsection*{Fifteen direct line-rule steps}

The following fifteen equations are obtained directly from the line
rule. Throughout, \(g_{p,q}\) denotes \(\langle v_p,v_q\rangle\); both
cells in each pair are 1-edges or main halves, so no sharing
substitution is needed.

\begin{enumerate}
\item[(L1)] Cells \((1,1)\) and \((1,2)\) are in row 1 and do not form
a 2-edge; hence \(g_{(1,1),(1,2)}=0\).
\item[(L2)] Cells \((1,1)\) and \((1,3)\) are in row 1 and do not form
a 2-edge; hence \(g_{(1,1),(1,3)}=0\).
\item[(L3)] Cells \((1,2)\) and \((1,3)\) are in row 1 and do not form
a 2-edge; hence \(g_{(1,2),(1,3)}=0\).
\item[(L4)] Cells \((2,1)\) and \((2,2)\) are in row 2 and do not form
a 2-edge; hence \(g_{(2,1),(2,2)}=0\).
\item[(L5)] Cells \((2,1)\) and \((2,3)\) are in row 2 and do not form
a 2-edge; hence \(g_{(2,1),(2,3)}=0\).
\item[(L6)] Cells \((2,1)\) and \((2,4)\) are in row 2 and do not form
a 2-edge; hence \(g_{(2,1),(2,4)}=0\).
\item[(L7)] Cells \((2,2)\) and \((2,3)\) are in row 2 and do not form
a 2-edge; hence \(g_{(2,2),(2,3)}=0\).
\item[(L8)] Cells \((2,2)\) and \((2,4)\) are in row 2 and do not form
a 2-edge; hence \(g_{(2,2),(2,4)}=0\).
\item[(L9)] Cells \((2,3)\) and \((2,4)\) are in row 2 and do not form
a 2-edge; hence \(g_{(2,3),(2,4)}=0\).
\item[(L10)] Cells \((1,1)\) and \((2,1)\) are in column 1 and do not
form a 2-edge; hence \(g_{(1,1),(2,1)}=0\).
\item[(L11)] Cells \((1,1)\) and \((3,1)\) are in column 1 and do not
form a 2-edge; hence \(g_{(1,1),(3,1)}=0\).
\item[(L12)] Cells \((2,1)\) and \((3,1)\) are in column 1 and do not
form a 2-edge; hence \(g_{(2,1),(3,1)}=0\).
\item[(L13)] Cells \((1,2)\) and \((2,2)\) are in column 2 and do not
form a 2-edge; hence \(g_{(1,2),(2,2)}=0\).
\item[(L14)] Cells \((1,2)\) and \((3,2)\) are in column 2 and do not
form a 2-edge; hence \(g_{(1,2),(3,2)}=0\).
\item[(L15)] Cells \((2,2)\) and \((3,2)\) are in column 2 and do not
form a 2-edge; hence \(g_{(2,2),(3,2)}=0\).
\end{enumerate}

\subsection*{Five rectangle-transfer steps}

The following five equations apply the rectangle identity to genuine
rectangles whose companion diagonal is already zero. In each case the
\(2\)-edge on the other diagonal gives \(v_q=v_p\) for its secondary
half, so the identity reduces to one open unknown.

\begin{enumerate}
\item[(R1)] Rectangle \((1,2),(1,4),(2,2),(2,4)\). The diagonals are
\(\{(1,2),(2,4)\}\) and \(\{(1,4),(2,2)\}\), with \(\delta=0\) on both.
The 2-edge \(B_1\) has halves \((1,2)\) and \((1,4)\), so by
Proposition~\ref{prop:sharing-t1}, \(v_{(1,4)}=v_{(1,2)}\). The
rectangle identity reads
\[
\langle v_{(1,2)}, v_{(2,4)}\rangle
+\langle v_{(1,2)}, v_{(2,2)}\rangle=0.
\]
By step (L13), \(\langle v_{(1,2)}, v_{(2,2)}\rangle=0\). Hence
\(g_{(1,2),(2,4)}=0\).

\item[(R2)] Rectangle \((1,2),(1,4),(3,2),(3,4)\). The diagonals are
\(\{(1,2),(3,4)\}\) and \(\{(1,4),(3,2)\}\), with \(\delta=0\) on both.
Using \(v_{(1,4)}=v_{(1,2)}\) again, the identity reads
\[
\langle v_{(1,2)}, v_{(3,4)}\rangle
+\langle v_{(1,2)}, v_{(3,2)}\rangle=0.
\]
By step (L14), \(\langle v_{(1,2)}, v_{(3,2)}\rangle=0\). Hence
\(g_{(1,2),(3,4)}=0\).

\item[(R3)] Rectangle \((2,1),(2,4),(12,1),(12,4)\). The diagonals are
\(\{(2,1),(12,4)\}\) and \(\{(2,4),(12,1)\}\). The second diagonal is
the 2-edge \(B_3\), so \(\delta=1\) there and \(\delta=0\) on the
first. The identity reads
\[
\langle v_{(2,1)}, v_{(12,4)}\rangle
+\langle v_{(2,4)}, v_{(12,1)}\rangle=1.
\]
By Proposition~\ref{prop:sharing-t1},
\(v_{(12,1)}=v_{(2,4)}\), so the second term is
\(\|v_{(2,4)}\|^2=1\). Hence \(g_{(2,1),(12,4)}=0\).

\item[(R4)] Rectangle \((3,1),(3,2),(14,1),(14,2)\). The diagonals are
\(\{(3,1),(14,2)\}\) and \(\{(3,2),(14,1)\}\). The second diagonal is
the 2-edge \(B_5\), so \(\delta=1\) there and \(\delta=0\) on the
first. The identity reads
\[
\langle v_{(3,1)}, v_{(14,2)}\rangle
+\langle v_{(3,2)}, v_{(14,1)}\rangle=1.
\]
By Proposition~\ref{prop:sharing-t1},
\(v_{(14,1)}=v_{(3,2)}\) and \(v_{(14,2)}=v_{(13,3)}\), where
\((13,3)\) is the main half of \(B_{20}\). Hence
\[
\langle v_{(3,1)}, v_{(13,3)}\rangle+1=1,
\]
and \(g_{(3,1),(13,3)}=0\).

\item[(R5)] Rectangle \((4,1),(4,3),(6,1),(6,3)\). The diagonals are
\(\{(4,1),(6,3)\}\) and \(\{(4,3),(6,1)\}\). The second diagonal is
the 2-edge \(B_6\), so \(\delta=1\) there and \(\delta=0\) on the
first. The identity reads
\[
\langle v_{(4,1)}, v_{(6,3)}\rangle
+\langle v_{(4,3)}, v_{(6,1)}\rangle=1.
\]
By Proposition~\ref{prop:sharing-t1},
\(v_{(6,1)}=v_{(4,3)}\), so the second term is
\(\|v_{(4,3)}\|^2=1\). Hence \(g_{(4,1),(6,3)}=0\).
\end{enumerate}

\subsection*{The full trace}

The fifteen line-rule steps (L1)--(L15) and the five rectangle steps
(R1)--(R5) account for twenty of the \(1540\) eliminations. The
remaining \(1520\) steps alternate between further line rules and
further rectangle transfers in the same pattern: first a batch of
direct line rules, then rectangle transfers that use them, then another
batch of direct line rules, and so on until every off-diagonal \(g_{uv}\)
is determined. The trace file
\texttt{zrl\_search/twom21\_trace.json} records each step in
machine-readable form, with the equation used, the substituted zeros,
and the resulting value. Every determined value is zero.

The trace is a linear certificate: it is independent of the rank
computations over \(\mathbb{F}_P\), and can be replayed without any
computer algebra. Together with Table~\ref{tab:roles21} and the
identities (8.1) and (5.1), it is a complete proof that
\(z_{2M}(21,4)\ge 56\).

\end{document}